\documentclass[conference]{IEEEtran}
\ifCLASSINFOpdf
   \usepackage[pdftex]{graphicx}
\else
\fi
\usepackage[cmex10]{amsmath}
\begin{document}
%
\title{A Context-Free Smart Grid Model Using Complex System Approach}

\author{\IEEEauthorblockN{Soufian Ben Amor}
\IEEEauthorblockA{
       University of Versailles SQY\\
       Versailles, France\\
              Email: soufian.benamor@uvsq.fr}
\and
\IEEEauthorblockN{Alain Bui}
\IEEEauthorblockA{University of Versailles SQY\\
       Versailles, France\\
       Email: alain.bui@uvsq.fr}
\and
\IEEEauthorblockN{Guillaume Gu\'erard}
\IEEEauthorblockA{University of Versailles SQY\\
       Versailles, France\\
       Email: guillaume.guerard@prism.uvsq.fr}}


%


\maketitle

\begin{abstract}

Energy and pollution are urging problems of the 21th century. By gradually changing the actual power grid system, smart
grid may evolve into different systems by means of size, elements and 
strategies, but its fundamental requirements and objectives will not
change such as optimizing production, transmission, and consumption. 
Studying the smart grid through modeling and simulation provides
us with valuable results which cannot be obtained in real world due to
time and cost related constraints.
 Moreover, due to the complexity of the smart grid, achieving global optimization is not an easy task. In this paper, we propose a complex system 
based approach to the smart grid modeling, accentuating on the optimization
by combining game theoretical and classical methods in different levels.
Thanks to this combination, the optimization can be achieved with flexibility
and scalability, while keeping its generality. 
\end{abstract}

\IEEEpeerreviewmaketitle

\section{Introduction}

Our society is electrically dependent. The electrical grid supply energy to households, businesses, and industries, but disturbances and blackouts are becoming common.  With the pressure from ever increasing energy demand
and climate change, finding new energy resources and enhancing energy 
efficiency have become priority of many nations in the 21st century.

The term smart grid is coined by Amin in 2005 \cite{amin2005toward}.
Smart grid is a type of electrical grid which attempts to predict 
and intelligently respond to the behavior and actions of all electric 
power users connected to it - suppliers, consumers and those that do 
both - in order to efficiently deliver reliable, economic, and sustainable 
electricity services.
Then, The expression ``Smart Grid'' has expanded into different dimensions: 
some see it as a numerical solution for downstream counter and 
mostly residential customers, while others believe that it is a global 
system vision that transcends the current structure of the energy 
market to generate economic, environmental, and social benefits for everyone.

Thus, Smart Grid is a fuzzy concept with various definitions in literature. However, Smart Grid could be defined according to the main requirements of an energy network. Smart Grid should integrate information and communication technologies to generate, transport, distribute, and consume energy more efficiently. In addition, the network should have mainly the following properties: self-healing, flexibility, predictability, interactivity, optimality, and safety \cite{guerard2012survey}. Moreover, the Smart Grid should improve reliability, reduce peak demand, and equalize energy consumption.

Research works are being conducted to attain the objectives, 
but many problems of modeling and coordination hamper advancements. 
However each model offers its own 
vision of the smart grid, putting aside theoretical and technological 
advancement of others. 
Cooperation between smart technologies and existing infrastructure is 
often neglected in scientific and industrial studies \cite{Molderink2009simulating}. 
In \cite{Borenstein2002dynamic}, authors argued that an electrical grid which 
allows the adjustments on both supply and demand will improve efficiency,
reduce costs on both sides and will be beneficial for the environment.

Taking into account all these internal and external features, the Smart Grid is defined as a complex system \cite{ahat2013smart, gao2012survey, guerard2012survey}. 
Contribution of our approach consists in
treating the smart grid as a complex system, locating the problems
at local as well as global levels, and solving them with coordinated
methods. In other words, through studying and analyzing
smart grid, we isolate homogeneous parts with similar behaviors or objectives, 
and apply classical optimization
algorithms at different levels with coordination. 
Thanks to combining those interdependent methods, our approach guarantees the 
flexibility in terms of
system size. Besides that, generality of our approach allows its 
applicability in different scenarios and models.

This paper is organized as following: in the next section, the concept of complex system is introduced and theoretical approaches in their modeling are discussed. In Section 3, we present the details of our global Smart Grid model based on the complex system approach, and in section 4, we present the research of a global consensus between supply and demand. We also discuss our perspectives and first results in section 5.

\section{Complex system approach}
A system which consists of large populations of connected
agents, or collections of interacting elements, is
said to be complex if there exists an emergent global
dynamics resulting from the actions of its parts rather
than being imposed by a central controller. That is a self-organizing
collective behavior difficult to anticipate from
the knowledge of local behavior \cite{boccara2004modeling}. Complex system
study embraces not only traditional disciplines of science,
but also engineering, management, and medicine \cite{turner2011symbiotic}.


The majority of studies on Smart Grids use a top-down approach. The Smart Grid is broken into basic issues: optimization, network structure and communication technologies and security \cite{fang2011smart}; and main components: users (consumers and producers), energy, controllers and data \cite{mas4}. Next, a global objective function is defined. The simulation returns the global solution, without taking into account local constraints. Applying the optimization method in complex
systems in a global manner is almost an impossible
task, if not impossible at all. Because, complex systems
are composed of heterogeneous parts; it is hard to find all the variables that matters; even if all variables are included, the complexity of the objective function will become beyond the computation power of computers \cite{kirkpatrick1983optimization}.

In addition, each organized activity shows a conflict between two fundamental requirements, allocation of resources in various tasks and coordination of these tasks to accomplish the global mission. While modeling complex systems, bottom-up analysis, also named systemic approach, gives a more complete, realistic and comprehensive vision \cite{macal2005tutorial}.
 
Smart grid can be qualified as a complex system \cite{guerard2012survey}, due to its heterogeneous actors, dynamic, complex interactions among them, and global behaviors such self-healing, or self-organizing characters. Based on this observation, we will analyze the smart grid methodologically in order to understand the mechanisms and internal components as well as the needs of every sub-components.

At first step, we should understand the system. An overview brings structural aspects, entities and objectives. All these elements are considered as agents. These one are not randomly distributed in the system, but according to patterns, and form distinct groups with their own arrangement. In the smart grid, three types of behavior are distinguished: consumer, producer, and transporter. The result is a hierarchy, or a food chain. 

After analyzing the characteristics of the system, we define the sub-components. A sub-component has a structure, objectives and specific entities; although quantities or position in the system are variable. As a separate system, it has its own dynamics. It is then possible to solve it with an appropriate optimization method. Consumers are generally located at the chain end, in network tree. The top producers are located in a mesh network, reinforcing the grid, and are linked to the consumers by linear chain.

The sub-components are interacting, then you should take into account the I-O data for each method.  The stability of the model depends on localized optimization. It is necessary to optimize each part of the chain as well as a whole to stabilize the system. If only consumers are optimized, all their devices will receive energy. If we optimize producer, they will produce the minimum at minimum cost. To prevent system crashes, the model must have a system of communication to reach a global consensus. Moreover, the system is subjected to external pressure. Feedback between sub-components are essential to maintain functionality, and to find local and global equilibrium.

These mechanisms work like homeostasis in natural sciences: a balance between the internal environment and the external environment. The Smart Grid, through its self-organization, self-healing and optimizing resources at any scale has a similar phenomenon. That is why the search for a global consensus is essential.

In summary, we analyze the system to determine sub-components. These have a system of communication and their own optimization methods, global criteria ensure balance in the system. Local optimization and global consensus constitute a decentralized optimization of our complex system guaranteeing individual and collective benefits.

\section{Modeling}
The problems of electrical networks have been known for long, and research as well as industrial works has been carried out to find effective and competitive solutions. Nevertheless the efforts are often concentrated on specific cases, and solutions are, too, specific without any room for evolution. Among the proposed solutions we can mention:
\begin{itemize}
\item Distributed generation/microgrids: since a centralized optimization is very costly in terms of time and memory, optimization should be done at all levels. The microgrids can change the centralized interface into a distributed interface, therefore optimization can be carried out in a distributed manner. Consequently calculation benefits in terms of time and memory are significant, while ensuring optimal at different scales.
\item Design of intelligent network (home automation): domotics or smart devices give real-time data and are controllable by the user or a smart meter. While optimizing local consumption, they optimize overall consumption as a result.
\item Energy storage device: the energy storage coupled with energy optimization from beginning to end, regulates consumption and clears consumption peaks.
\item Reduction of Transmission and Distribution T\&D network losses by automated distribution: One of the strong points of our model is the distribution optimization by local and global algorithms which reduce the loss of congestion or routing errors.
\item Intelligent control of price: when the network be- comes intelligent, it is necessary that the consumer prices may also change in order to follow the new consumer behavior.

\end{itemize}

Many theories, in order to optimize smart grid, come from complex system analysis.
Sub-components, optimization methods and the necessary theories to obtain an optimal consensus, are defined in these articles \cite{ahat2013smart,guerard2012survey,guerard2012}.

\subsection{Global objective function}

 The overall mathematical problem is similar to a knapsack problem \cite{martello1990knapsack,xiong2011smart}. The objective function is under multiple temporal, spatial and physical constraints: such as the granularity of study, local optimization (routing, distribution, consumption), as well as variable consumption and production time.

\begin{figure}[!ht]
\centering
\includegraphics [width=2.8in] {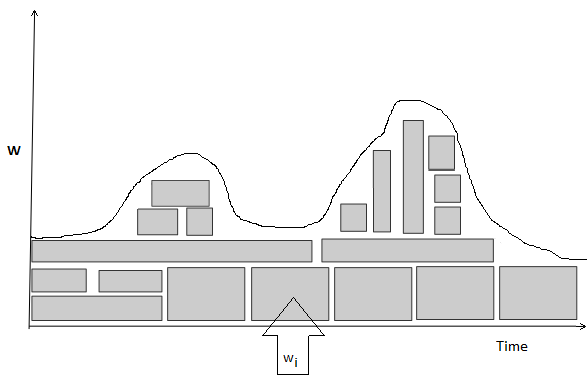}
\caption{Knapsack problem for Smart Grid.}
\label{fig1}
\end{figure}

The general 0-1 knapsack problem applied to smart grids is:
$$
\left\{
\begin{array}{l}
\mbox{ maximize } \;\sum_{i=1}^{n}x_i u_i,\\
\mbox{ subject to } \sum_{i=1}^{n}x_i w_i \leq W
\end{array}
\right.
$$ 
where $x_i = 1$ if the devices $i$ is taken, else $0$; $u_i$ represents utility of the device $i$ and $w_i$ its consumption in $Watt*hour$; $W$ is the total energy produced in $Watt*hour$; there are $n$ devices. Moreover, several quadratic or linear constraints due to the complex system is added (routing, minimal values, cost, etc.).

This problem, see Figure \ref{fig1}, is too hard to be solved at very large scale - millions of items, in few minutes. Criteria or constraints should be satisfied throughout algorithms and process.
Is the decomposition of the global problem and assigning it to various computers connected to a single network is equivalent to resolving the global solution on a single machine? In distributed algorithms, all machines have the same role. We notice that all levels have the same overall goal, but only use specific algorithms. In other words, our process is similar to distributed algorithms.

\subsection{A three layered grid} 

The Smart Grid has three sub-components having a structure, dynamics and distinct behavior: the transmission and distribution network (T\&D), the microgrid and the local level.

\begin{figure}[!ht]
\centering
\includegraphics [width=2.2in]{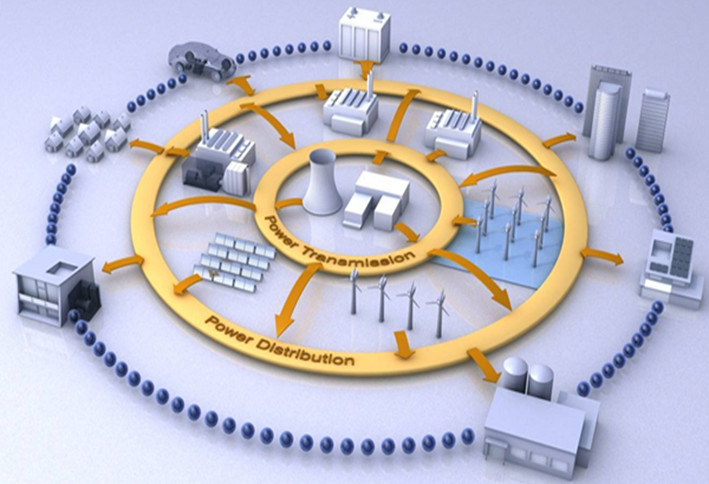}
\caption{Smart Grid sub-components (from PowerMatrix, Siemens).}
\label{fig2}
\end{figure}

The first level is the only full connected one, forming a single group, represented by the center, see Figure \ref{fig2}. The Transmission and Distribution network (T\&D) must deliver energy from producers to points of consumption. Energy flows on electric cables with various criteria and technical constraints limiting the amount of energy that can circulate on each. The algorithm at this level should be able to limit the effects of congestion n due to the widespread use of a few lines, while limiting the cost of routing energy. Production and consumption must match as better as possible, in order to achieve this, we must deliver most of energy while satisfying most of consumers

The second level is the link between consumption and energy production, represented by the second ring, see Figure \ref{fig2}. The microgrid is a broader view of local consumers, it is a structure representing an eco-district bounded by the upstream substation. Its role is to distribute energy from substation to consumers. For this, it books an amount of energy from the T\&D network.

The outer ring represents local levels, see Figure \ref{fig2}. Local level models consumers, i.e. a group of consuming devices, local renewables or electric vehicle requesting or providing a measurable amount of energy.  These isolated structures, like residences or factories support the consumption of energy, i.e. the distribution of energy among appliances under its responsibility. The energy distribution is a dynamic programming resolution of a knapsack problem: objects are devices, the weight is the energy received at the local level, and the utility is a function of demand-side management strategies

\subsection{Iteration process}
\begin{figure}[!ht]
\centering
\includegraphics [width=3.4in]{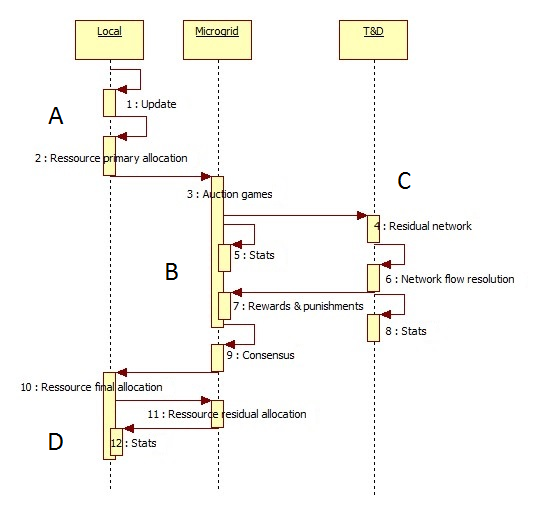}
\caption{Sequential Scheme.}
\label{fig4}
\end{figure}

An iteration occurs every five minutes. Once data are updated, the process is decomposed into four sequences, see Figure \ref{fig4}.

Sequence A: to design intelligent aspect of the device, a priority is assigned to these dynamic entities, and for calculating a consumption value. Indeed, we use a local knapsack problem, solved by dynamic programming after data normalization, for finding a primary first optimal resource allocation. The knapsack is solved by dynamic programming, its complexity is $O(n*C_{0})$ where $n$ is the number of devices and $C_{0}$ is the energy received. It is reduced to $O(n)$ if prognostics are correct. Due to the number of devices and the size of the bag, the optimal solution on a 0-1 knapsack by dynamic programming is obtained instantly with normalized weight.

Sequence B: this sequence aims to book an amount of energy from producers to consumers using an auction. There are two way to book energy: a consensus between consumers and production, i.e. a game where microgrids and energy flows are players; and a bid system with feedback. The problem of the first one is the complexity of the problem, impossible to resolve in few time. Second way have the advantage of time, but don't guarantee the global optimum. During an auction, it is likely that energy required does not correspond to any new consumption strategy, i.e. local level consume as much with that or without. In addition, it is possible to search the nearest consumer at lowest cost. The number of possible strategies is infinite, we must look at the impact of each of them on the microgrid and on the final decision.

\begin{figure}[!ht]
\centering
\includegraphics [width=2.5in]{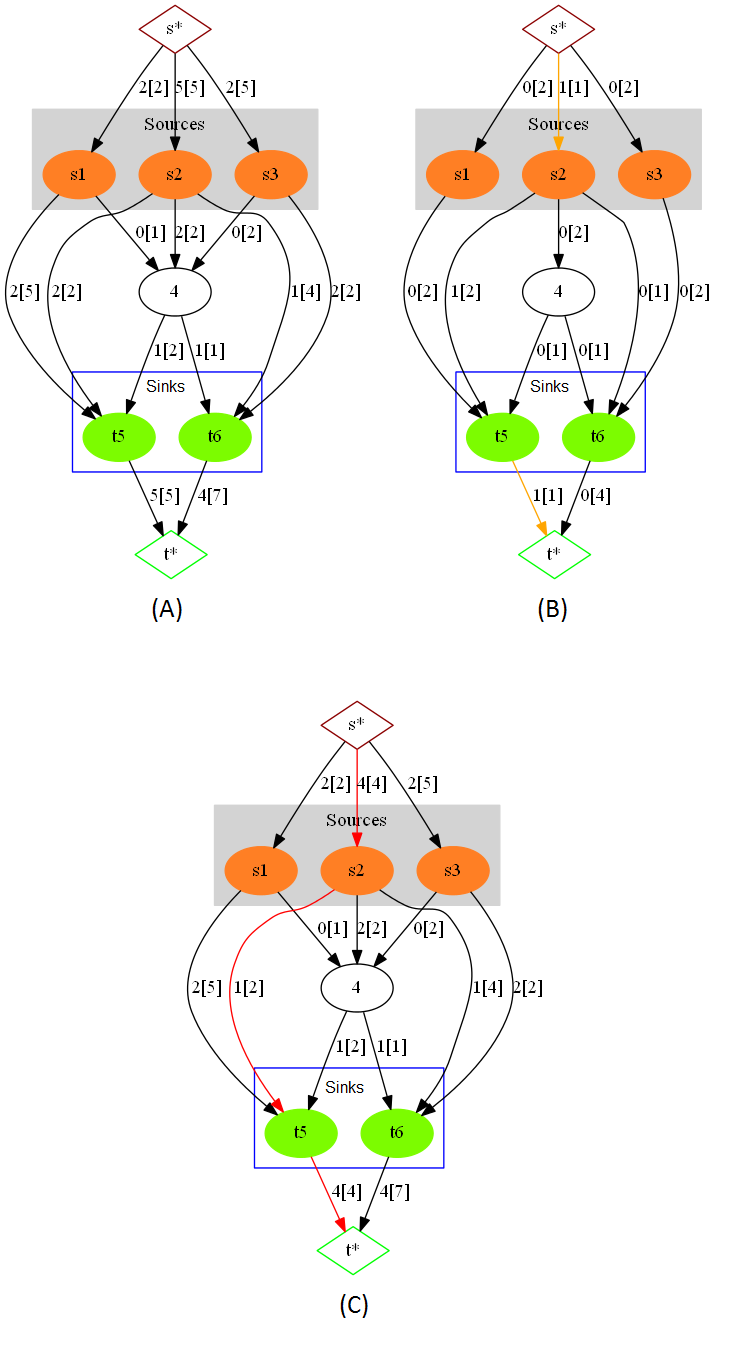}
\caption{Updating of routing.}
\label{rourou}
\end{figure}

Sequence C: About the problem of routing, nodal rule or Kirchhoff's circuit specifies that at any node in a circuit, the sum of currents flowing into that node is equal to the sum of the currents flowing out of that node. An electrical circuit is equal to a graph in which a junction is a node, and physical connection corresponds to an edge. Routing problem is equivalent to the known max flow problem. Gale's theorem shows the existence of a solution in a network of offers and requests \cite{gale1957theorem}. The flow of the previous iteration is maximum by Ford-Fulkerson, recalculate the entire flow is not necessary. The residual graph removes excess flows between two updates, optimizing the computation time of the optimal flow, see Figure \ref{rourou}. It is also possible to resolve a maximum flow with minimal cost and minimal flow on edges by Busacker and Gowen, but this algorithm can't use previous and update graph solutions. If supply and demand do not match, algorithm analyzes the bottlenecks by performing Ford-Fulkerson on two schemes: infinite production, infinite consumption. These data are also used to calculate prognostics. The complexity is $O(A*f)$ where $A$ is the number of edges and $f$ is the maximum flow. The residual graph reduces $f$ to the sum of local differences in production and consumption.

Sequence D: energy is distributed by knapsack problem, according to the auctions. The unconsumed energy is redistributed among non-used devices at upper scale. The device's priorities are updated according to the result of the final distribution. At worst, the complexity is $O(n*K)$, with $n$ the number of devices and $K$ the energy received.

\subsection{Global direction}

Studies conducted by Barabasi \cite{barabasi2003linked} and Watts \cite{watts2004six} raised four general principles in distributed adaptive systems, more generally in complex systems:
\begin{enumerate}
\item Global information is encoded as statistics and dynamic patterns in the components of the system.
\item Chance and probabilities are essential.
\item The system performs a parallel search of opportunities.
\item The system has a continuous interaction \cite{mitchell2006complex}.
\end{enumerate}

Principles of electricity generation and distribution are well known. Synchronization of the system is recognized that each station and each piece of equipment runs on the same clock, which is crucial for its proper functioning. Cascading failures related desynchronization can lead to massive power outages. In smart grid, a technical control automates the management of energy. Real-time data must be converted into information quickly enough so that errors are diagnosed in time, corrective actions are identified and executed dynamically in the network, and feedback loops provide measures to ensure that the per- formed actions and production are consistent \cite{anderson2011adaptive}.

\begin{figure}[!ht]
\centering
\includegraphics [width=3.1in]{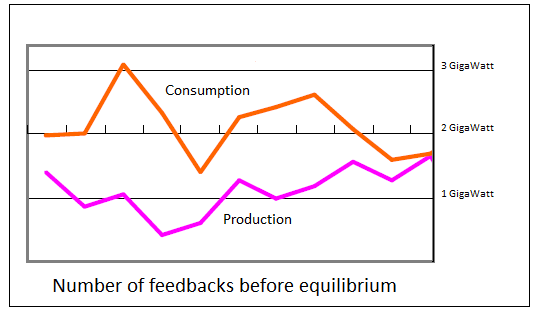}
\caption{Supply and demand's Consensus.}
\label{fig5}
\end{figure}

An arbitrary configuration generates a random pattern, without prognostics. As a result of the first auction, gap between consumption and production occurs. Feedback adjusts supply and demand, as renewable energies or electric vehicle management for example. It is easier to vary demand over supply. So, supply and demand tend to the same value, see Figure \ref{fig5}. There are many steady states, economic and optimal criteria are considered in the final solution, both seeking minimal cost.

Global regulation must be done both at consumer and producer level, taking into account the difficulties of routing. To increase the effectiveness of this method, it is assumed that the front part of the infrastructure is home automation and every device can be controlled separately by the user and regulation algorithms. Regulation is an overview of the Smart Grid in order to smooth the production curve.
We can discern three types of regulation:
\begin{enumerate}
\item
Mathematical regulations: mathematical tools are introduced to smooth the consumption curve (standard, planning according to the derivative, gradient and barycenter, etc.).
\item
Regulation by self-stabilization: the criteria for regulation of the curve are done at any point of the smart grid. Some technologies are already in effect, such as dynamic pricing systems or consumer subscriptions.
\item
Hybrid regulation: This type of regulation is based on both mathematical and self-stabilizing approach. Its main advantage is to minimize the risks associated with either method.
\end{enumerate}
The model is actually based on a mathematical regulation.

\section{Equilibrium between supply and demand}
\subsection{Demand-side management}

\begin{figure}[!ht]
\centering
\includegraphics [width=3.4in]{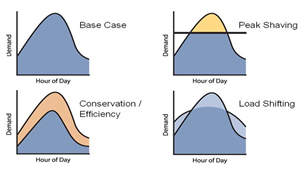}
\caption{Load shifting strategies.}
\label{dsm}
\end{figure}

\begin{figure}[!ht]
\centering
\includegraphics [width=3.2in]{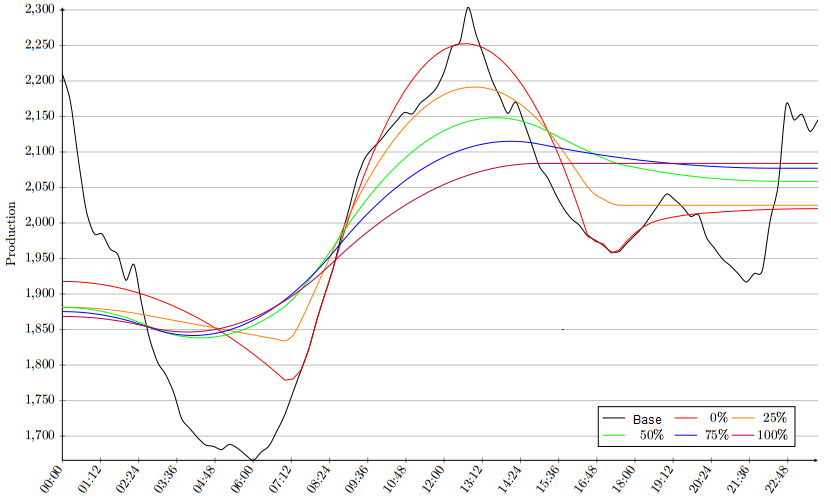}
\caption{Load shifting based on electric vehicles \cite{hermans2012}.}
\label{load}
\end{figure}

In order to predict consumption, Smart Grid will allow customers to make informed decisions about their energy consumption, adjusting both the timing and quantity of their electricity use. This ability to control usage is called demand-side management (DSM).
In the literature, DSM programs have two goals: demand-response programs for energy efficiency; and load shifting which schedules the production and consumption over the long term, see Figures \ref{dsm} and \ref{load}.

Initially, energy conservation programs encourage customers to give some energy use in return for saving money, such as turning up the thermostat a few degrees in summer to reduce air conditioning.  Additional gains in energy efficiency are possible through technologies that can provide targeted education or real time verification of costumer demand reduction. But the consumer behavior is too sporadic and cannot be implemented in the model by variables. However, it is possible to describe the desired effects at the microgrid through strategies. So, the auction will be based on a multitude of strategies for each consumer.

Demand response programs and load shifting transfer customer load during periods of high demand to off-peak periods. Shifting daily peak demand flattens the load curve, allowing more electricity to be provided by less expensive base load generation. 

DSM programs have existed across the globe since the 1970s and is an active field of research \cite{zachhuber2008simulating}. California utilities have used such programs to hold per-capita energy consumption nearly constant over the past 30 years.  DSM programs incorporate some or all of the following six levers \cite{davito2010smart}:
\begin{enumerate}
\item Rates or utility tariffs. Not yet implemented, we currently work with economists to define to price energy (depending on producers), and the impact of DSM on consumers and producers prices.
\item Incentives. Consumers’ participation in demand-side management programs is simulated by strategies. The impact on priority and the price isn't fix and depend on cases.
\item Access to information. Local level, microgrid and T\&D level have access to information and algorithms allow energy management. 
\item Utility controls, simulated by priority and algorithms process.
\item Education and marketing.
\item Customer insight and verification, see subsection \ref{eco}.
\end{enumerate}

\subsection{Strategies based on utility}
To explain the role of strategies, a simplistic case will be served as an example.
Let five houses (H1 to H5) are in the same microgrid, see Table \ref{tab1}. The following table provides devices in every home with their variables: in order consumption, priority of operation ($0$ meaning already in use) and the value for knapsack. Values of device $i$, noted $u_i$ in each house are calculated as follow : let $w_{max}$ the greatest consumption in the house, $p_{max}$ the greatest priority in the house; for each device $i$ in the house, $u_i = (w_{max} *  p_{max}) - (w_i *  p_i) + w_i$. The table also present the forecast (Fc.) and the minimal energy required (devices in bold).

\begin{table}[h]
\begin{center}
\caption{Consumption in the microgrid.}
\begin{tabular}{|r|r|r|r|r|r|}
\hline
 & H1 & H2 & H3 & H4 & H5\\
\hline
Dev. & \textbf{1/0/81} & \textbf{1/0/16} & \textbf{1/0/0} & \textbf{1/0/33} & \textbf{1/0/0}\\
 & \textbf{1/1/80} & \textbf{1/0/16} & \textbf{1/0/0} & \textbf{1/1/32} & \textbf{3/0/0}\\
 & \textbf{3/0/83} & \textbf{2/1/15} & \textbf{10/0/0} & \textbf{3/0/35} &  \\
 & 5/2/75 & \textbf{3/0/18} & & 3/2/29 &  \\
 & 20/4/20 & 4/3/7 &  & \textbf{4/1/32} &  \\
 &  & 5/3/5 &  & 8/4/8 &  \\
Fc. & 4 & 6 & 12 & 8 & 6\\
Min. & 5 & 7 & 12 & 9 & 4\\
\hline
\end{tabular}
\label{tab1}
\end{center}
\end{table}

Information are send to microgrid. Let $l$ is the utility of the strategy for the consumers, and $r$ for the producer. Different DSM strategies are defined for each house:
\begin{enumerate}
\item Basic consumption: consumption of all devices, all combination possible. In the example, combination are based on the priority of the devices. Utilities are calculated as following: $l=\sum_{i=1}^{n}\dfrac{u_i * w_i}{p_i}$ for each device $i$ in this strategy; $r=\sum_{i=1}^{n}(\dfrac{u_i}{p_i} -\alpha)*w_i$ with $\alpha$ the average utility for an unity of consumption.
\item Peak shaving. Priority have an exponential impact on the utility of consumer and producer.
\item Conservation. All utilities depend on priority except those that can provide energy.
\item Load shifting. Utilities depend on the average time of consumption of all devices and the total amount of energy needed.
\item Over-production. Priority of batteries are reduced, in order to reload them.
\item Over-consumption. If possible, batteries give energy, and domotic reduces its consumption.
\end{enumerate}
We present the result of simplified basic consumption strategies in the Table \ref{tab2}. Each strategy are based on a priority level, i.e. all devices equal or less than the priority are taken into account. Results are presented as follows: house value noted $r$/distribution value noted $l$). At the end of the table, the final consumption (after three feedback) is shown. It is calculated with the strategy with the max $r+l$.

\begin{table}[h]
\begin{center}
\caption{Strategies of the microgrid.}
\begin{tabular}{|r|r|r|r|r|r|}
\hline
Priority & H1 & H2 & H3 & H4 & H5\\
\hline
0 & 330/194 & 77/27 & \textbf{Done} & 138/47 & \textbf{Done}\\
1 & 410/240 & \textbf{107/37} &  & 298/56.5 & \\
2 & \textbf{620/257} & none &  & \textbf{341/32.5} &  \\
3 & none & 125/-36 & & none &  \\
4 & 720/-322 & none &  & 357/-131 &  \\
5 & none & none &  & none &  \\
Final & 10 & 7 & 12 & 12 & 4\\
\hline
\end{tabular}
\label{tab2}
\end{center}
\end{table}

Currently strategies are unilateral, only the behavior of the local level is taken into account. To avoid many feedback, the market economy is studied in order to find in few games the final result. It is not intended to promote the producer but to plan the routing during a cooperative game.

\subsection{Economic goals}\label{eco}
Like any investment decision in technology, the benefits need to exceed the costs. Making the value 
case for smart grid investments is complicated by at least two characteristics. First, the smart grid assets 
 contribute to more than one value stream. Making a 
value determination for smart grid investments 
usually requires the recognition and accounting 
of the benefits from multiple value streams to 
offset the investment costs in technology 
deployment. Second, several of these value 
streams can be difficult to quantify financially. 
Reliability is traditionally something that is set 
by regulation and best practice and implemented 
as a necessary cost of 
providing electricity \cite{apec2006apec}.
Determining the value of decreasing environmental impact and ensuring the health and well-being of the 
populace are examples of other areas where benefits from smart grid investment are hard to capture in 
simple equations \cite{pratt2010smart}.

\begin{figure}[!ht]
\centering
\includegraphics [width=3.2in]{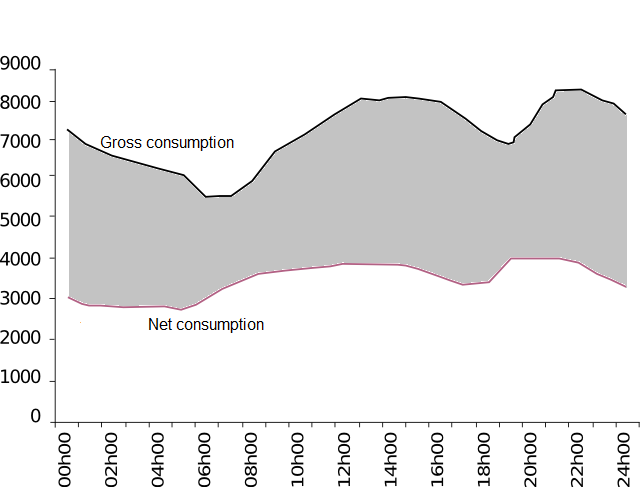}
\caption{Difference between normal consumption and consumption with management and renewable energies.}
\label{gfdgdf}
\end{figure}

Nevertheless, it is possible to measure the level of profits curves of Smart Grid consumption compared to standard curves. Indeed, consumption curves are known and used for many years to schedule daily production. 
At the local level, gross consumption, i.e. consumption without the aid of any control technology or renewable energy used locally, is compared to net consumption, i.e. consumption in the Smart Grid, see Figure \ref{gfdgdf}. The cost of all renewable energy, local or plants, is approximate. Profits are the difference between the cost of the gross consumption minus the net consumption, and the cost of the technologies used.

\begin{figure}[!ht]
\centering
\includegraphics [width=3.2in]{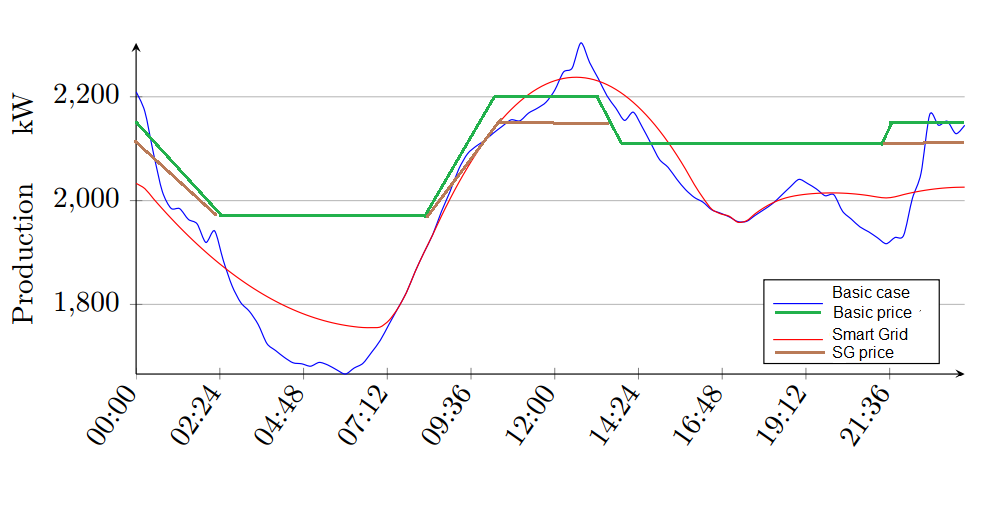}
\caption{Energy prices before and after DSM.}
\label{gezgez}
\end{figure}

Similarly, the cost of fuel plants or other plant used during peak consumption is known. On a global scale, the average cost of energetic output current is approximate, in the basic case, and in a Smart Grid, see Figure \ref{gezgez}.

\section{Experimental results and Discussions}

\begin{figure}[!ht]
\centering
\includegraphics [width=3.4in]{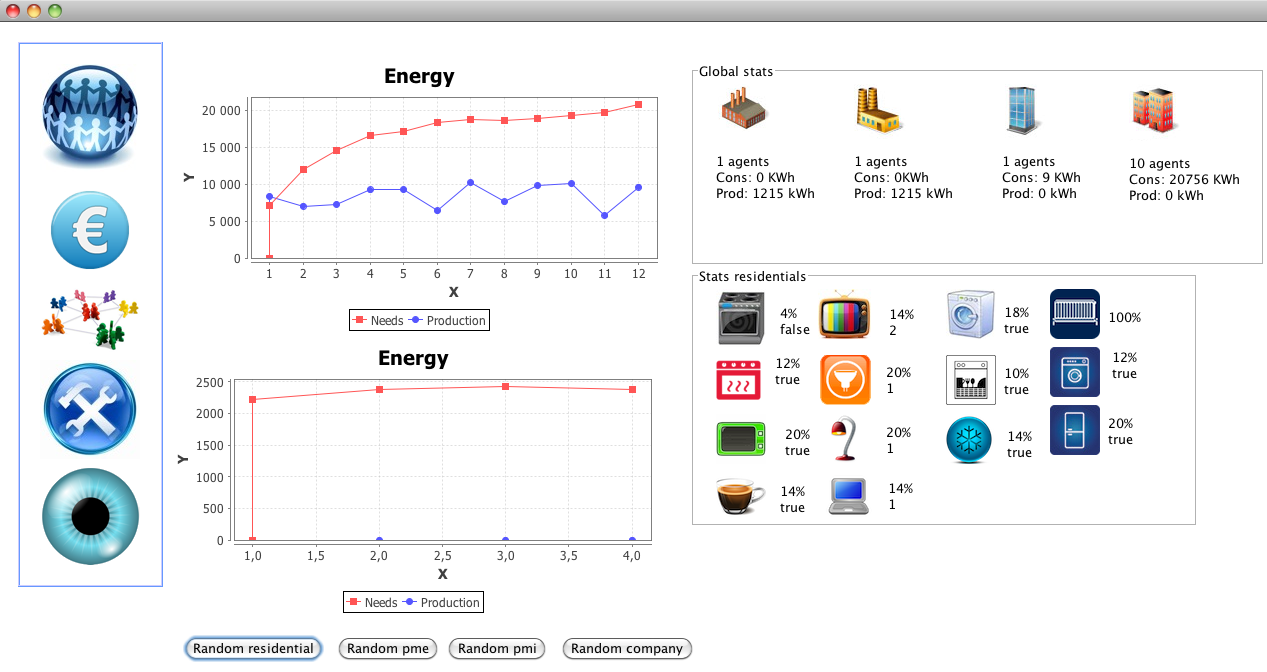}
\caption{Overview of the local production and consumption.}
\label{fig8}
\end{figure}

The model was implemented using a multi-agent simulation paradigm \cite{ferscha2009efficiency, jennings1998roadmap, petermann2012complex} (Figure \ref{fig8}). Each agent is either a consumer, a producer, both or an energy's transporter. Agents have specific behaviors induced by their class and act depending on previous algorithms.
To validate the model, instances at local and global scale have been made. Agents present like engine consumption or energy plants' production are implemented by French national production companies data and energy distribution data (EDF and RTE) \cite{BBPe12c}. 

In the first tests, consumption and production tend towards equilibrium. Local and renewable energies are privileged to maximize their profitability. The model limits the losses of the distance of consumption, and uses the least amount of fossil energy. It works at any scale and any agents under the condition that there exists a feasible solution. First results are based on an arbitrary utilities that does not take into account the economic aspect. The economic study will do a forthcoming publication.

About the mathematical global problem, the knapsack for the overall system, we provide two very simple and thus highly practical algorithms that found the global solution of convex maximization quadratic problems.
Algorithms are based respectively on inner and outer approximation like ball or cuboid following by a local search and the standard cutting plane technique adapted to our problem. The papers about these methods are not yet published.

The goal of the Smart Grid's model is to reduce difference between the mathematical result and the model's result. Learning strategies are being developed so that the Smart Grid is the best possible configuration without human intervention.

The proposed model works for randomized or parameterized Smart Grids, we actually work in the Positive Energy 2.0 project led by ALSTOM Energy Management and various companies such as Bouygues or Renault to validate the model on real projects.

\section{Conclusion}
As smart grid can be qualified as a complex system, classical optimization methods cannot be applied directly, due to the computational complexity in terms of time and memory.
More generally, we also demonstrated how to solve optimization problems in complex systems. While applying optimization algorithms directly in complex systems is nearly impossible, we should analyze the system and divide them into sub-systems with defined characteristics, then we should apply specific algorithms and coordinate them using multi-agent simulation in order to achieve global optimization.

A general context-free model of a smart grid is being developed, which integrates those algorithms. This model does not replace the current model nor provides an ideal model, but presents an improvement of the Energy Grid. Preliminary tests have validate our approach.
Data mining and learning strategy are studied in order to limit the number of simulation and variable's modifications before to obtain an optimized configuration.

%
\bibliographystyle{abbrv}
\bibliography{single}  

\end{document}